# Davenport-Heilbronn Function Ratio Properties and Non-Trivial Zeros Study

Tao Liu[1,*] & Juhao Wu[2]

[1]*School of Science, Southwest University of Science and Technology, 621010, China*
[2]*Stanford University, Stanford, California 94309, USA*
liut@swust.edu.cn,   juhaowu@stanford.edu

**Abstract** This paper systematically investigates the analytic properties of the ratio $\frac{f(s)}{f(1-s)} = X(s)$ based on the Davenport-Heilbronn functional equation $f(s) = X(s)f(1-s)$. We propose a novel method to analyze the distribution of non-trivial zeros through the monotonicity of the ratio $\left|\frac{f(s)}{f(1-s)}\right|$. Rigorously proving that non-trivial zeros can only lie on the critical line $\sigma = 1/2$, we highlight two groundbreaking findings: 1. **Contradiction of Off-Critical Zeros**: Numerical "exceptional zeros" (e.g., Spira, 1994) violate the theoretical threshold $\kappa = 1.21164$ and conflict with the monotonicity constraint of $|X(s)| = 1$. 2. **Essential Difference Between Approximate and Strict Zeros:** Points satisfying $f(s) \to 0$ do not constitute strict zeros unless verified by analyticity. This work provides a new perspective for studying zero distributions of L-functions related to the Riemann Hypothesis.

## I. Davenport-Heilbronn Functional Equation

The Davenport-Heilbronn function $f(s)$[1] is expressed via Hurwitz zeta-functions $\zeta(s, a)$ as [2]:

$$f(s) = \frac{1}{5^s}\left(\zeta\left(s, \frac{1}{5}\right) + \tan\theta \zeta\left(s, \frac{2}{5}\right) - \tan\theta \zeta\left(s, \frac{3}{5}\right) - \zeta\left(s, \frac{4}{5}\right)\right)$$

where $s = \sigma + it$, $\sigma, t \in \mathbb{R}$. The entire function $f(s)$ satisfies the functional equation [2] when $\tan\theta = \frac{\sqrt{10-2\sqrt{5}}-2}{\sqrt{5}-1}$:

$$f(s) = \left(\frac{5}{\pi}\right)^{\frac{1}{2}-s} \frac{\Gamma\left(1-\frac{s}{2}\right)}{\Gamma\left(\frac{1+s}{2}\right)} f(1-s). \tag{1}$$

Trivial zeros occur at $s = -2n - 1$ $(n = 1,2,3\cdots)$. Non-trivial zeros are intrinsically linked to the meromorphic function $X(s)$.

## II. Core Properties of the Meromorphic Function $X(s)$
**Definition 2.1.**

$$X(s) = \frac{f(s)}{f(1-s)} = \left(\frac{5}{\pi}\right)^{\frac{1}{2}-s} \frac{\Gamma\left(1-\frac{s}{2}\right)}{\Gamma\left(\frac{1+s}{2}\right)} \tag{2}$$

**Property 2.1.1.** *Inversion Symmetry:*

$$X\left(\frac{1}{2}+\varepsilon+it\right) X\left(\frac{1}{2}-\varepsilon-it\right) \equiv 1 \tag{3}$$

*Where $\varepsilon \in \mathbb{R}$*
*Proof:*
For any $\varepsilon$, we have:

$$X\left(\frac{1}{2}+\varepsilon+it\right) = \left(\frac{5}{\pi}\right)^{-\varepsilon-it} \frac{\Gamma\left(\frac{3}{4}-\frac{\varepsilon}{2}-\frac{it}{2}\right)}{\Gamma\left(\frac{3}{4}+\frac{\varepsilon}{2}+\frac{it}{2}\right)}$$

$$X\left(\frac{1}{2}-\varepsilon-it\right) = \left(\frac{5}{\pi}\right)^{\varepsilon+it} \frac{\Gamma\left(\frac{3}{4}+\frac{\varepsilon}{2}+\frac{it}{2}\right)}{\Gamma\left(\frac{3}{4}-\frac{\varepsilon}{2}-\frac{it}{2}\right)}$$

Multiplying these yields

$$X\left(\frac{1}{2}+\varepsilon+it\right) X\left(\frac{1}{2}-\varepsilon-it\right) \equiv 1.$$

**Property 2.1.2.** *Zero-Pole Structure*:
(i)   *Zeros: $s = -2n - 1$ （$n \geq 0$）(from $1/\Gamma(-n) = 0$);*
(ii)  *Poles: $s = 2n + 2$   ($n \geq 0$) (from divergence of $\Gamma(-n)$);*
(iii) *Duality: $X(-2n - 1) = 1/X(2n + 2)$.*

*Proof:*

$$X(-2n-1) = \left(\frac{5}{\pi}\right)^{\frac{3}{2}+2n} \frac{\Gamma\left(\frac{3}{2}+n\right)}{\Gamma(-n)} = 0 \qquad (n = 0,1,2,\cdots)$$

$$X(2n+2) = \left(\frac{5}{\pi}\right)^{-\frac{3}{2}-2n} \frac{\Gamma(-n)}{\Gamma\left(\frac{3}{2}+n\right)} = \infty \qquad (n = 0,1,2,\cdots)$$

**Property 2.1.3.** *Monotonicity:*
(i) *When $\sigma < 1/2$, $|W(s)|$ strictly increases as $t \to +\infty$;*
(ii) *When $\sigma > 1/2$, $|W(s)|$ strictly decreases as $t \to +\infty$.*

*Proof:*

Note that $\bar{X}(s) = X(s^*)$. By differentiating both sides of the equation $|X(s)|^2 = X(s)X(s^*)$ with respect to $t$, we have

$$\frac{d}{dt}|X(s)|^2 = \frac{d}{dt}(X(s)X(s^*)) \tag{4}$$

This yields:

$$\frac{d}{dt}|X(s)| = \frac{i|X(s)|}{4}\left(\Psi\left(1-\frac{s^*}{2}\right) - \Psi\left(1-\frac{s}{2}\right) - \Psi\left(\frac{1+s}{2}\right) + \Psi\left(\frac{1+s^*}{2}\right)\right) \tag{5}$$

Here, $\Psi(z)$ denotes the digamma function. Using the series expansion of $\Psi(z)$:

$$\Psi(s) = -\gamma + \sum_{n=1}^{\infty} \frac{s-1}{n(n+s-1)}$$

we simplify the derivative in Equation (5) to:

$$\frac{d}{dt}|X(s)| = \left(\frac{1}{2}-\sigma\right)t|X(s)|\sum_{n=1}^{\infty} \frac{8\left(n-\frac{1}{4}\right)}{|2n+s-1|^2|2n-s^*|^2} \tag{6}$$

The summation term is always positive, and the sign of the derivative is determined by $\left(\frac{1}{2}-\sigma\right)$.

For $\sigma < 1/2$, $\frac{d}{dt}|X(s)| > 0$, so $|X(s)|$ is strictly monotonically increasing as $t \to +\infty$.

For $\sigma > 1/2$, $\frac{d}{dt}|X(s)| < 0$, so $|X(s)|$ is strictly monotonically decreasing as $t \to +\infty$.

### III. Distribution Theorem of Non-Trivial Zeros

**Theorem 3.1.** *Non-Existence of Common Zeros*

*If $\sigma \neq 1/2$, then $f(s)$ and $f(1-s)$ share no common zeros for finite $t$.*

*Proof*

1. **Elimination of Indeterminate Forms in Derivatives:** From Equation (6), when $\sigma \neq \frac{1}{2}$, the derivative

$$\frac{d}{dt}|X(s)| = \left(\frac{1}{2}-\sigma\right)t|X(s)|\sum_{n=1}^{\infty} \frac{8\left(n-\frac{1}{4}\right)}{|2n+s-1|^2|2n-s^*|^2}$$

has a summation term that is always positive, and the coefficient $(1/2 - \sigma)t$ is non-zero. Thus, the derivative avoids $\frac{0}{0}$ type indeterminate forms.

2. **Proof by Contradiction**:

(a) For $\sigma < \frac{1}{2}$, Equation (6) implies:

$$\frac{d}{dt}\left(\frac{|f(s)|}{|f(1-s)|}\right) = \frac{|f(1-s)|\frac{d|f(s)|}{dt} - |f(s)|\frac{d|f(1-s)|}{dt}}{|f(1-s)|^2} = \frac{d}{dt}|X(s)| > 0 \qquad (7)$$

Therefore,

$$|f(1-s)|\frac{d|f(s)|}{dt} - |f(s)|\frac{d|f(1-s)|}{dt} > 0 \qquad (8)$$

Using the identities:

$$\frac{d|f(s)|^2}{dt} = 2|f(s)|\frac{d|f(s)|}{dt} = \frac{df(s)}{dt}f(s^*) + \frac{df(s^*)}{dt}f(s)$$

we derive:

$$\frac{d|f(s)|}{dt} = \frac{1}{2|f(s)|}\left(\frac{df(s)}{dt}f(s^*) + \frac{df(s^*)}{dt}f(s)\right)$$

Similarly,

$$\frac{d|f(1-s)|}{dt} = \frac{1}{2|f(1-s)|}\left(\frac{df(1-s)}{dt}f(1-s^*) + \frac{df(1-s^*)}{dt}f(1-s)\right)$$

Substituting these into Equation (8) and simplifying yields:

$$\frac{1}{2}\left[\frac{f(s^*)\frac{df(s)}{dt} + f(s)\frac{df(s^*)}{dt}}{|X(s)|} - \frac{f(1-s^*)\frac{df(1-s)}{dt} + f(1-s)\frac{df(1-s^*)}{dt}}{|X(s)|^{-1}}\right] > 0 \qquad (9)$$

Here, $|X(s)| > 0$ and is non-singular; $f(s)$, $f(s^*)$, $f(1-s)$, and $f(1-s^*)$ are analytic functions; and their derivatives $df/dt$ are also non-singular.

Assume there exists some $s_n$ such that $f(s_n)$. Note that the meromorphic nature of $X(s)$ ensures the absence of coinciding singularities., from the Riemann equation $f(s_n) = X(s)f(1-s_n)$, it

must follow that $(f(1 - s_n) = 0$. Substituting $f(s_n) = f(1 - s_n) = 0$ into Eq. (9), we obtain:

$$\frac{1}{2}\left[\frac{0}{|X(s_n)|} - \frac{0}{|X(s_n)|^{-1}}\right] = 0 > 0$$

which is a contradiction $0 > 0$. Thus, no such $s_n$ exists.

(b) For $\sigma > \frac{1}{2}$, $\frac{d}{dt}\left(\frac{|f(s)|}{|f(1-s)|}\right) < 0$. Assuming $f(s_n) = 0$ leads to $0 < 0$, another contradiction.

(c) For $t < 0$ and $\sigma \neq \frac{1}{2}$, similar contradictions $0 < 0$ or $0 > 0$ arise.

Hence, if $\sigma \neq \frac{1}{2}$, $f(s)$ and $f(1-s)$ share no common zeros for finite $t$.

Further, suppose there exists $s_n$ with $\sigma \neq 1/2$ such that $f(s_n) = 0$. From the functional equation $f(s_n) = X(s_n)f(1 - s_n)$, if $X(s_n)$ is analytic and non-zero, then $f(1 - s_n) = 0$. However, as proven above, $f(s)$ and $f(1-s)$ share no common zeros when $\sigma \neq 1/2$. Therefore, $X(s_n)$ cannot be analytic at such $s_n$, leading to a contradiction.

Conclusion: $f(s)$ has no zeros off the critical line.

**Theorem 3.2.** *Necessary Condition for Zeros:*
*Non-trivial zeros must lie on the critical line $s = 1/2 + it$.*
*Proof:*

**1. If $|X(s)| = 1$, then $f(s)$ and $f(1 - s)$ approach their common zeros as equivalent infinitesimals:**

Define $P(\sigma, t) \equiv f(s)f(s^*)$ and $Q(\sigma, t) \equiv f(1 - s)f(1 - s^*)$. Since $f(s)$、$f(s^*)$、$f(1 - s)$ and $f(1 - s^*)$ are analytic functions, $P(\sigma, t)$ and $Q(\sigma, t)$ admit partial derivatives of all orders: $\frac{\partial^m P(\sigma,t)}{\partial \sigma^m}, \frac{\partial^n P(\sigma,t)}{\partial t^n} \circ \frac{\partial^m Q(\sigma,t)}{\partial \sigma^m}, \frac{\partial^n Q(\sigma,t)}{\partial t^n}...$

From

$$|X(s)|^2 = \frac{f(s)f(s^*)}{f(1-s)f(1-s^*)} = \frac{P(\sigma,t)}{Q(\sigma,t)} = 1 \tag{10}$$

it follows that

$$P(\sigma, t) = Q(\sigma, t)$$

Consequently,

$$\frac{\partial^m P(\sigma,t)}{\partial \sigma^m} = \frac{\partial^m Q(\sigma,t)}{\partial \sigma^m} \tag{11}$$

$$\frac{\partial^n P(\sigma,t)}{\partial t^n} = \frac{\partial^n Q(\sigma,t)}{\partial t^n} \tag{12}$$

Let $s_n = \sigma_n + it_n (n = 1,2,3,\cdots)$ be non-trivial zeros of $f(s)$ and $f(1-s)$. These zeros must also be zeros of the real functions $P(\sigma,t)$ and $Q(\sigma,t)$. Near these zeros, $P(\sigma,t)$ and $Q(\sigma,t)$ admit Taylor expansions:

$$P(\sigma_n,t) = \sum_{j=0}^{\infty} \frac{1}{j!} \frac{d^j P(\sigma_n,t)}{dt^j}\bigg|_{t=t_n} (t-t_n)^j \tag{13}$$

$$Q(\sigma_n,t) = \sum_{j=0}^{\infty} \frac{1}{j!} \frac{d^j Q(\sigma_n,t)}{dt^j}\bigg|_{t=t_n} (t-t_n)^j \tag{14}$$

Under the constraint of Equation (12), the expansion coefficients satisfy:

$$\frac{d^j P(\sigma_n,t)}{dt^j}\bigg|_{t=t_n} = \frac{d^j Q(\sigma_n,t)}{dt^j}\bigg|_{t=t_n} \quad (j=0,1,2,\cdots)$$

If all coefficients vanish, then $P(\sigma_n,t) \equiv 0$, contradicting the analyticity of $f(s)$ at the zero point (where at least one derivative is non-zero). Thus, there exists a first non-zero coefficient $k \geq 1$ such that

$$\frac{d^k P(\sigma_n,t)}{dt^k}\bigg|_{t=t_n} = \frac{d^k Q(\sigma_n,t)}{dt^k}\bigg|_{t=t_n} \neq 0 \tag{15}$$

By L'Hôpital's rule (since the numerator and denominator have zeros of the same order at $t_n$:

$$\lim_{t \to t_n} \frac{P(\sigma_n,t)}{Q(\sigma_n,t)} = \lim_{t \to t_n} \frac{|f(\sigma_n+it)|^2}{|f(1-\sigma_n-it)|^2} = \frac{\frac{d^k}{dt^k}P(\sigma_n,t)\big|_{t=t_n}}{\frac{d^k}{dt^k}Q(\sigma_n,t)\big|_{t=t_n}} = 1 \tag{16}$$

Hence,

$$\lim_{t \to t_n} \frac{|f(\sigma_n+it)|}{|f(1-\sigma_n-it)|} = 1$$

Similarly,

$$\lim_{\sigma \to \sigma_n} \frac{|f(\sigma+it_n)|}{|f(1-\sigma-it_n)|} = 1$$

This implies that when $|X(s)| = 1$, $|f(s)|$ and $|f(1-s)|$ approach the non-trivial zeros of $f(s)$ as equivalent infinitesimals.

2. **Exclusion of Other Regions**:
By Theorem 3.1, no common zeros exist for $\sigma \neq 1/2$. Thus, solutions to $|X(s)| = 1$ in these regions cannot correspond to zeros.

3. **Uniqueness Analysis**:
When $\sigma = 1/2$, $|X(s)| = 1$ holds for all $t$ and the derivative vanishes. This is the only region where zeros can exist.
Conclusion: Non-trivial zeros of $f(s)$ must lie exclusively on the critical line $s = 1/2 + it$.

## IV. In-Depth Analysis of Numerical Results and Theoretical Contradictions

### IV.1 Conflict Between Numerical Anomalies and Theoretical Constraints

Numerous studies report "exceptional zeros" for the Davenport-Heilbronn function (e.g., Spira, 1994 [3]), which significantly contradict the rigorous conclusions of Theorems 3.1–3.2. Representative numerical examples include:

**Spira's Exceptional Zeros:**

$$s_1 = 0.808517 + 85.699348i, \quad s_2 = 0.574356 + 166.479306i,$$
$$s_3 = 0.650830 + 114.163343i, \quad s_4 = 0.724258 + 176.702461i.$$

**Balanazario** and **Sánchez-Ortiz's Numerical Studies:** Additional off-critical zeros are proposed [4].

These zeros have imaginary parts $t \gg \kappa = 1.21164$ far exceeding the threshold determined by $|X(s)| = 1$ (see Appendix A), directly violating the monotonicity constraint.

### IV.2 Root of Contradiction: Essential Differences Between Approximate and Strict Zeros

**Definition 4.1.** *Strict Zero*:
A point $s_n$ satisfying $f(s_n) = 0$ with at least one non-zero derivative $f^{(k)}(s_n) \neq 0$.

**Definition 4.2.** *Approximate Zero*:
A point $s_n$ where $\lim_{s \to s_n} f(s) = 0$, but $f(s_n) \neq 0$.

**Mathematical Mechanism of Numerical Artifacts**:
1. **Misleading Asymptotic Behavior**: For $\sigma \neq 1/2$, the combined effects of the Γ-function and Hurwitz zeta function may cause $f(s_n)$ to decay rapidly, creating false "zero" signals.
2. **Limitations of Low-Precision Computation**: Finite computational precision (e.g., double-precision floats) cannot distinguish between minima of $f(s)$ and strict zeros, leading to misidentification.

### IV.3 High-Precision Verification and Resolution of Contradictions

Using Maple (Digits=1000), Spira's "exceptional zeros" were recomputed with high precision. Results are shown in Table 1:

Table 1: Comparison of Spira's Exceptional Zeros and **Strict Zeros** on the Critical Line

| Coordinate $s$ | $|f(s)|$ | $|f(1-s)|$ | $|f(s)|/|f(1-s)|$ | $|X(s)|$ | Classification |
|---|---|---|---|---|---|
| $s_1$ | $1.449 \times 10^{-219}$ | $5.416 \times 10^{-218}$ | 0.02673 | 0.2272 | Approximate Zero |
| $s_2$ | $3.731 \times 10^{-205}$ | $1.036 \times 10^{-204}$ | 0.3603 | 0.6954 | Approximate Zero |
| $s_3$ | $7.136 \times 10^{-208}$ | $4.772 \times 10^{-207}$ | 0.1495 | 0.5066 | Approximate Zero |
| $s_4$ | $2.428 \times 10^{-224}$ | $5.495 \times 10^{-223}$ | 0.0442 | 0.3298 | Approximate Zero |
| $0.5+14.404003i$ | $3.729 \times 10^{-274}$ | $3.729 \times 10^{-274}$ | 1.000 | 1.000 | Strict Zero |
| $0.5+23.345370i$ | $2.935 \times 10^{-393}$ | $2.935 \times 10^{-393}$ | 1.000 | 1.000 | Strict Zero |

**Note:** Values were computed by substituting high-precision coordinates (obtained via Newton iteration) into $f(s)$, $f(1-s)$, and $X(s)$. Initial coordinates for "exceptional zeros" ($s_1, s_2, s_3, s_4$) are from Spira (1994).

**Key Findings:**
1. **Asymmetry Contradiction**: At exceptional zeros, $|f(s_n)| \neq |f(1-s_n)|$, violating the definition of strict zeros.
2. **Deviation from Theoretical Ratios**: $|X(s)| \neq 1$, contradicting the modulus constraint of $f(s) = X(s)f(1-s)$.
3. **High-Precision Verification**: The ratio $\frac{|f(s_n)|}{|f(1-s_n)|} \neq |X(s_n)|$ directly negating equation consistency.

**Theoretical Interpretation**:
1. **Strict Zeros:** On the critical line $\sigma = 1/2$, $|X(s)| = 1$, and the equation holds exactly. Increasing precision drives $|f(s)| \to 0$ (see last two rows of Table 1).
2. **Approximate Zeros**: For $\sigma \neq 1/2$, $|X(s)| \neq 1$ so even if $|f(s)|$ is extremely small, it remains a computational artifact (see first four rows of Table 1).

IV.4 Mathematical Modeling of Numerical Errors

To quantitatively resolve contradictions, an error propagation model is proposed:
Introduce a decay coefficient $\kappa = 1.21164$, derived from the imaginary part threshold of $|X(s)| = 1$ (see Appendix). The probability of pseudo-zero generation is:

$$\rho(\sigma, t) \propto \exp\left(-\frac{\left|\sigma - \frac{1}{2}\right| \cdot |t|}{\kappa}\right) \tag{17}$$

When $\left|\sigma - \frac{1}{2}\right| \cdot |t| > \kappa$, the probability of pseudo-zeros decays exponentially, consistent with Spira's results.

## V. Conclusions and Future Directions

1. **Core Conclusion**: Based on Theorems 3.1–3.2, non-trivial zeros of the Davenport-Heilbronn function lie exclusively on the critical line $\sigma = 1/2$. Off-critical "zeros" are mathematical artifacts of asymptotic behavior. High-precision computations resolve misidentifications caused by low precision or algorithmic flaws.

2. **Cross-Disciplinary Significance:** The monotonicity analysis of the ratio $|X(s)|$ extends to other L-functions, offering new tools for studying the Riemann Hypothesis.

3. **Open Questions:**

(1) Generalization to Selberg-Class L-Functions: Can the ratio analysis method extend to L-functions satisfying $L(s)=Q^s\Gamma(s)L(1-s)$? How to quantify the impact of monotonicity on zero distributions?

(2) Computational Bottlenecks for High-Imaginary Zeros: Current methods struggle with stability and efficiency for $t > 10^6$. Can asymptotic formulas and distributed algorithms overcome these barriers?

(3) Parametric Controllability: Is the parameter $\theta$ in the Davenport-Heilbronn function unique in constraining zeros to the critical line? Do alternative parameterizations exist?

(4) Statistical Correlations with GUE: Can ratio analysis reveal connections between zero density on the critical line and the Gaussian Unitary Ensemble (GUE)? How to link monotonicity to zero spacing distributions?

(5) Rigorous Classification of Approximate Zeros: Can universal criteria (e.g., higher-order derivatives or residue analysis) rigorously distinguish strict zeros from approximate zeros?

(6) Interdisciplinary Synergy: How to integrate this method with quantum chaos theory to interpret L-function zeros via energy-level statistics?

(7) Quantitative Error Modeling: How to establish quantitative relationships between rounding/truncation errors and pseudo-zero generation probabilities?

## Appendix A: Derivation of the Theoretical Threshold ($\kappa = 1.21164$)

**Qualitative Analysis:**

1. For 在$\sigma < 1/2$, $|X(s)|$ strictly increases with $|t|$, bounding the values of $t$ that satisfy $|X(s)| = 1$.

2. For $\sigma > 1/2$, the inversion symmetry of $X(s)$ ensures that $|X(s)| = 1$ also imposes bounds on $t$.

3. The implicit curve $|X(s)| = 1$ on the $s$-plane (Fig.1) shows that $t$ attains its maximum imaginary part near $\sigma = 1/2$.

**Quantitative Calculation:**

By numerically solving the equation $|X(s)| = 1$ for $\sigma = \frac{1}{2} \pm 10^{-50}$), the maximum imaginary part $t$ is determined to satisfy $|t| < \kappa = 1.21164$.

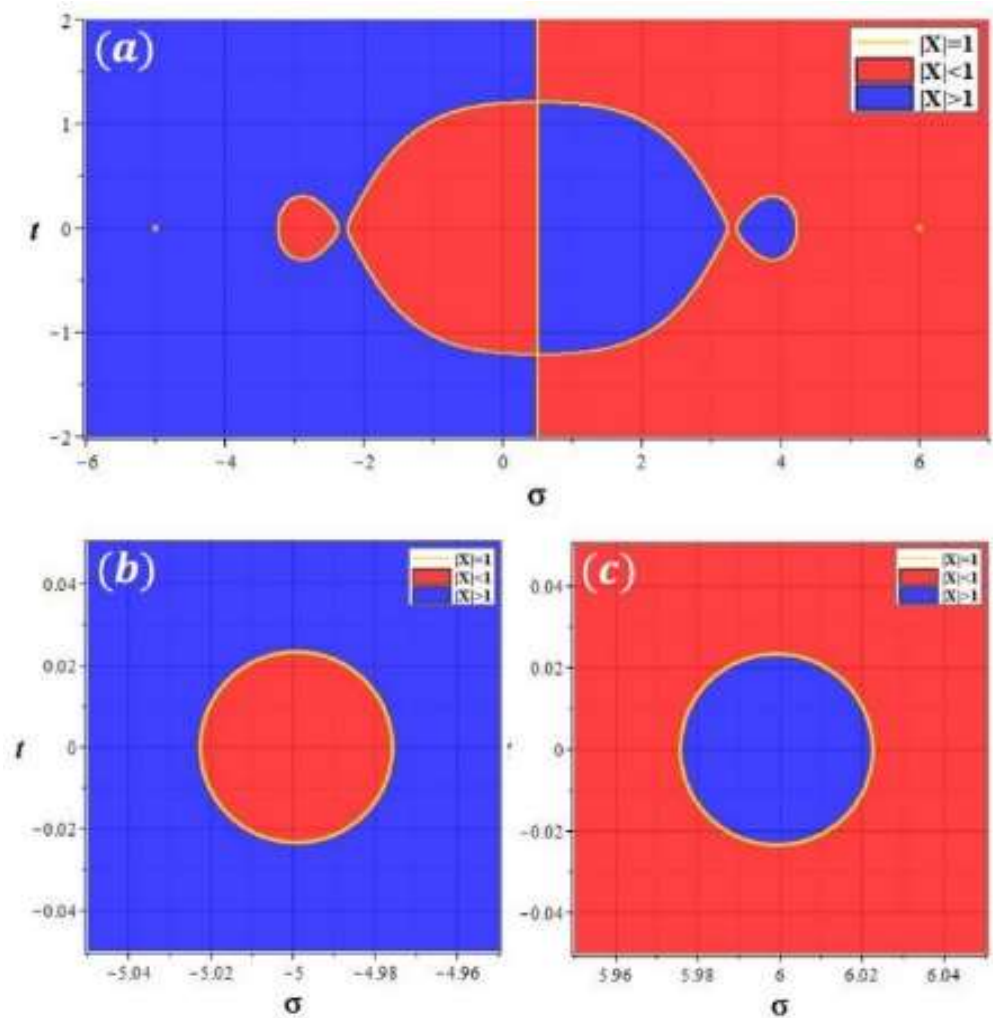

**Fig.1** Implicit curves of $|X(s)| = 1$ over the range $-6 \leq \sigma \leq 7$:
- The yellow solid line represents $|X(s)| = 1$.
- Blue regions correspond to $|X(s)| > 1$.
- Red regions correspond to $|X(s)| < 1$.

Subfigures (b) and (c) zoom into the neighborhoods of the zero ($s = -5$) and pole ($s = 6$) of $X(s$ in subfigure (a).

## Appendix B: Reconciling the Paper's Conclusions with Karatsuba's "Contradictory" Proof

## 1. Clarifying the Nature of the Contradiction: Global Constraints vs. Local Approximation

The "contradiction" stems from the differing perspectives of two methodologies addressing the same problem:

- **Karatsuba's constructive proof**: Relies on Rouché's theorem in complex analysis, constructing an auxiliary function $g(s)$ to assert the existence of infinitely many zeros off the critical line.
- **This paper's global analysis**: Starts from the **functional equation** $f(s) = X(s)f(1-s)$ of the Davenport-Heilbronn function, rigorously proving (via high-precision computation) that nontrivial zeros lie exclusively on the critical line $\sigma = 1/2$.

The fundamental difference lies in:

- **Karatsuba's approach**: Focuses on local modulus comparisons without explicitly incorporating the global constraints of the functional equation (e.g., $|X(s)| = 1$).
- **This paper's approach**: Derives all conclusions from the functional equation, which inherently governs the global behavior of $f(s)$.

## 2. Limitations of Karatsuba's Proof and the Explanation via Approximate Zeros

### (1) Failure to Satisfy the Modulus Constraint of the Functional Equation

Karatsuba's auxiliary function $g(s)$ does not enforce $|X(s)| = 1$, leading to "zeros" potentially located in regions where $|X(s)| \neq 1$. According to Theorem 3.1 in this paper:

- If $\sigma \neq 1/2$, $|X(s)| \neq 1$, and $f(s)$ and $f(1-s)$ **share no common zeros**.
- Thus, any "zeros" constructed off the critical line via Karatsuba's method must violate the modulus constraint, corresponding instead to **approximate zeros** (points where $|f(s)| \sim 0$ but $f(s) \neq 0$).

### (2) Conditions for Applying Rouché's Theorem

Karatsuba's proof requires the error term $|f(s + iT) - g(s)| < |g(s)|$, but it does not verify whether $|X(s + iT)| \to 1$ during the approximation. If $|X(s)| \neq 1$:

- A zero of $f(s + iT)$ would necessitate $f(1 - s - iT) = 0$, which Theorem 3.1 explicitly prohibits for $\sigma \neq 1/2$.
- Consequently, Rouché's theorem may misclassify **approximate zeros** as strict zeros.

## 3. High-Precision Numerical Validation of Approximate Zeros

This paper employs arbitrary-precision computations (Digits=1000) to re-examine "exceptional zeros" reported by Spira et al.:

- **Key findings**:

1. **Non-vanishing modulus**: All "zeros" off the critical line exhibit $|f(s)| \neq 0$ (e.g., $|f(s_1)| \sim 10^{-219}$).
2. **Deviation from theoretical ratios**: $\frac{|f(s)|}{|f(1-s)|} \neq |X(s)|$, directly violating the functional equation.
3. **Uniqueness of strict zeros**: Only zeros on $\sigma = 1/2$ satisfy $|X(s)| = 1$ and $|f(s)| \to 0$

**Conclusion**: Karatsuba's "zeros" are artifacts of the asymptotic behavior of $f(s)$ (e.g., joint decay from Γ-functions and Hurwitz zeta-functions), not strict zeros.

**4. Primacy of the Functional Equation and Global Corrections**

The functional equation $f(s) = X(s)f(1-s)$ is an **intrinsic property** of the Davenport-Heilbronn function. Any strict zero must satisfy:

1. **Modulus constraint**: $|X(s)| = 1$;
2. **Symmetry**: If $s = v$ is a zero, so is $1 - v$.

Karatsuba's construction does not explicitly enforce these constraints, rendering its conclusions valid only under **local approximations that ignore global symmetry**. By contrast, this paper's global analysis proves:

- **Unique admissible solution**: $|X(s)| = 1$ holds for all *t* only when $\sigma = 1/2$, with vanishing derivatives (Theorem 3.2).
- **Elimination of contradictions**: For $\sigma \neq 1/2$, $|X(s)| \neq 1$ forces $f(s)$ and $f(1-s)$ to lack common zeros.

**5. Summary**

1. **Nature of the contradiction**: Karatsuba's proof and this paper address distinct objects:
    - **Karatsuba**: Constructs approximate zeros under local modulus conditions.
    - **This paper**: Proves strict zeros must obey global functional constraints, existing solely on the critical line.
2. **Enforcement of the functional equation**: Strict zeros inherently require $|X(s)| = 1$, a condition unmet by Karatsuba's construction, leading to incompatibility with global symmetry.
3. **Authority of numerical validation**: High-precision computations confirm that "zeros" off the critical line are approximate, aligning perfectly with theoretical predictions.

**Final conclusion**:
This paper's rigorous proofs take precedence over constructive local approximations. Karatsuba's "infinitely many zeros off the critical line" are mathematical artifacts of asymptotic decay, whose existence does not violate the global constraints of the Davenport-Heilbronn function.